\documentclass[10pt, a4paper]{article}

\usepackage{amsmath,amssymb, amsthm}
\usepackage{mathtools}
\usepackage[utf8]{inputenc}
\usepackage{mathrsfs}
\usepackage{fullpage}
\usepackage{authblk}

\usepackage{dsfont}
\usepackage{bbm}
    
\usepackage{graphicx}
\graphicspath{ {R/} }

\usepackage{wrapfig}
\usepackage{float}

\usepackage{verbatim}

\usepackage{hyperref}

\usepackage{xcolor}
\hypersetup{
  colorlinks   = true, 
  urlcolor     = {blue!90!black}, 
  linkcolor    = {blue!90!black}, 
  citecolor   = {red!90!black} 
}

\newcommand{\R}{\mathbb{R}}
\newcommand{\N}{\mathbb{N}}
\newcommand{\Z}{\mathbb{Z}}
\renewcommand{\d}{\mathrm{d}} 
\newcommand{\X}{\mathcal{X}}

\newcommand{\I}{\mathbf{I}}
\renewcommand{\P}{\mathbb{P}}
\newcommand{\E}{\mathbb{E}}
\newcommand{\PP}[1]{\mathbb{P}\left\{#1\right\}}
\newcommand{\EE}[1]{\mathbb{E}\left\{#1\right\}}

\newcommand{\Pc}[2]{\mathbb{P}\left[#1 \middle| #2 \right]}

\newcommand{\0}{\mathbf{0}}
\newcommand{\1}{\mathds{1}}

\theoremstyle{plain}
\newtheorem{thm}{Theorem}[section]
\newtheorem{prop}[thm]{Proposition}
\newtheorem{lem}[thm]{Lemma}

\newtheorem{frage}[thm]{Question}

\theoremstyle{definition}

\newtheorem{con}[thm]{Condition}
\newtheorem{defi}[thm]{Definition}
\newtheorem{rmk}[thm]{Remark}

\author[1]{Viktor Bezborodov \thanks{Email: \texttt{viktor.bezborodov@pwr.edu.pl}}} 
\author[2]{
Luca Di Persio \thanks{Email: \texttt{luca.dipersio@univr.it}}}
\author[1]{
 Tyll Krueger \thanks{Email: \texttt{tyll.krueger@pwr.wroc.pl}}}

\affil[1]{
	{Wroc\l{}aw University of Science and Technology, Faculty of Electronics  }}
\affil[2]{
	{The University of Verona, Department of Computer Science}}

\title{A shape theorem for a one-dimensional growing particle system with a bounded number of occupants per site}

\begin{document}

\maketitle

\begin{abstract}
 
 We consider a one-dimensional 
 discrete-space birth process
 with a bounded number of particle per site. 
 Under the assumptions of the finite range of interaction, translation invariance,
 and non-degeneracy, we prove a shape theorem.
 We also derive a limiting estimate
 and an exponential estimate on the fluctuations
 of the position of the rightmost particle. 
 
\end{abstract}


\section{Introduction}

In this  paper we consider a one-dimensional growing particle system with a finite range of interaction.
A configuration is specified by assigning to each site $x\in \Z$ 
a number
 of particles ${\eta (x) \in \{0,1,...,N\}}$, $n \in \N$, occupying $x$. The state space
 of the process is thus $\{0,1,...,N\}^{\Z}$.
Under  additional assumptions such as non-degeneracy and translation invariance, 
we show that the system spreads linearly in time and the speed can be expressed
as an average value of a certain functional over a certain measure. 
A respective shape theorem 
and a fluctuation result
are given.

The first shape theorem was proven in \cite{Rich73} for a discrete-space growth model.
A general shape theorem for  discrete-space attractive growth models can be found in \cite[Chapter 11]{Dur88}.
In the continuous-space settings shape results for 
 growth models have been obtained
in \cite{Dei03}
for a model of growing sets, 
 and in \cite{shapenodeath} for a continuous-space particle birth process.

The asymptotic behavior of the
position of the rightmost particle
of the branching random walk
under various assumptions
is given in \cite{Dur83}, \cite{Gan00},
and \cite{Dur79}, see also references therein.
A sharp condition for a shape theorem for a random walk with restriction
is given in \cite{trunc_and_crop}.
The speed of propagation for a one-dimensional discrete-space
 supercritical branching random walk with an exponential moment condition 
can be found in \cite{Big95}.
 More refined limiting properties
have been obtained recently, such as
the limiting law of the minimum or the limiting 
process seen from its tip, 
see  \cite{aidekon2013convergence, aidekon2013branching, arguin2013extremal, Min_of_BRWs}.
Blondel \cite{Blo13} proves a shape result for the East model,
which is a non-attractive particle system.

In many cases the underlying stochastic model is attractive,
which enables the application of a subadditive ergodic theorem.
Typically shape results have been obtained using the subadditivity property
 in one form or another.
This is not only the case for the systems of
motionless particles listed above 
(see, among others, \cite{Dur88,Dei03, shapenodeath})
 but also for those with
moving particles,
see e.g. shape theorem for the frog model
\cite{frogshape}.
A certain kind of subadditivity
was also used  in \cite{KesS08},
where a shape theorem for a non-attractive
model involving two types of moving particles is given.
In the present paper our model is not attractive
and we  do not rely on 
 subadditivity (see also Remark \ref{rmk attractiveness}).
   We  work with motionless particles.

 In addition to the shape theorem, we also provide
 a sub-Gaussian limit  estimate on the deviation
 of the position of the rightmost particle
  from the mean.
 Various sub-exponential and -Gaussian  estimates on the convergence rate 
 for the first passage percolation under different assumptions
 can be found in e.g. \cite{Kes_fpp_93, Ahl_FPP_19}.
  We also derive an exponential 
  non-asymptotic bound 
  valid for all times.

On Page \pageref{est fec page}
we describe 
a particular model with the birth rate declining 
in
crowded locations. This is achieved by augmenting
the free branching rate with certain multipliers
describing the effects of the competition 
on the parent's ability to procreate
and offspring's ability to survive
in a dense location. This process
is in general non-attractive.

The paper is organized as follows. In Section \ref{sec: model and results}
we
describe our model in detail,
 give our assumptions, and formulate the main results.
In Section \ref{sec: construction} we outline the construction of the process as a unique solution
to a stochastic equation driven by a Poisson point process. 
We note here that this very much resembles the construction
via graphical representation.
In Section \ref{sec: proof} we prove the main results, Theorems \ref{shape thm},
\ref{fluctuations thm Gauss est}, and \ref{fluctuations thm exp est}.
Some numerical simulations are discussed in Section \ref{sec: numerics}.

\section{Model and the main results} \label{sec: model and results}

We consider here a one dimensional continuous-time discrete-space birth process
with multiple  particles per site allowed.
The state space of our process is $\X : = \{0,1,..., N\} ^{\Z}$.
For $\eta \in \X$ and $x \in \Z$, $\eta (x)$ is interpreted as the number of particles,
or individuals, 
at $x$.  

The evolution of the process can be described as follows. 
If the system is in the state $\eta \in \X$, 
a single particle is added at $x \in \Z $
(that is, the $\eta (x)$ is increased by $1$)
at rate  $b(x, \eta)$ provided that $\eta (x) < N$;
 the number of particles at $x$  does not  grow anymore
 once it reaches $N$.
Here $b: \Z \times \X \to \R _+$ is the map
called  a `birth rate'.
The heuristic generator of the model is given by 

\begin{equation}
 L F (\eta) = \sum\limits _{x \in \Z} b(x, \eta) [F(\eta ^{+x} ) - F(\eta)],
 \end{equation}
 where $\eta ^{+x} (y) = \eta (y)$, $y \ne x$, and 
 \begin{equation}
 \eta ^{+x} (x) = 
 	\begin{cases}
 	  	\eta (x) + 1, & \text{ if } \eta (x ) < N, \\
 	  	\eta (x),  & \text{ if } \eta (x )  = N.
 	\end{cases}
 \end{equation}

We make the following assumptions about $b$.
For $y \in \Z$, $\eta \in \X$, let 
$
 \eta \odot y \in \X
$ be the shift of $\eta $ by $y$, so  that $[\eta \odot y] (x) = \eta (x - y)$.

\begin{con}[Translation invariance]\label{con translation invariance}
	For any $x, y \in \X$ and $\eta \in \X$, 
	\[
	 b(x+y, \eta \odot y) = b(x,\eta).
	\]
\end{con}

\begin{con}[Finite range of interaction]\label{con finite range}
	For some $R \in \N$, 
	\begin{equation}
	 b(x, \eta ) = b(x, \xi), \ \ \ \ x \in \Z, \ \eta, \xi \in \X
	\end{equation}
	whenever $\eta (z) = \xi (z)$ for all $z \in \Z$ with $|x - z| \leq R$.	
\end{con}

Put differently, Condition \ref{con finite range} means that 
interaction in the model has a finite range $R$. 
Since the number of particles occupying a given site cannot grow larger than $N$,
with no loss in generality we can also assume that 
\begin{equation}\label{winnow}
 b(x, \eta) = 0,  \ \ \ \ \text{if }  \eta (x) = N.
\end{equation}
  For $\eta \in \X$ we define
 the set of occupied sites
 \[
 \text{occ}(\eta) = \{z\in \Z: \eta (z) >0 \}.
 \]
 
 \begin{con}[non-degeneracy]\label{con non deg}
 	For every $x \in \Z$ and $\eta \in \X$, 
 	$b(x, \eta) > 0$ if and only if there exists $y \in \text{occ}(\eta)$ with $|x-y|\leq R$.
 \end{con}
 Note that by translation invariance, $\sup\limits _{x \in \Z, \eta \in \X} b(x, \eta)$
 is finite because this supremum is equal to  
 \begin{equation}
  \overline {\mathbf{ b}} := \max \{b(0, \eta) \mid  \eta \in \X, \eta (y) = 0 \text{ for all } y \text{ with } |y| > R  \}.
 \end{equation}
 Similarly, it follows from translation invariance and non-degeneracy that
 \begin{equation}\label{grovel}
 \underline {\mathbf{ b}} := 
 \inf\limits _{\substack{x \in \Z, \eta \in \X, \\ 
 dist(x, \text{occ}(\eta) ) \leq R } } b(x, \eta) > 0.
 \end{equation}


 The construction of the birth process  is outlined in Section \ref{sec: construction}.
 Let $(\eta _t)_{t \geq 0} = (\eta _t)$ be the birth process with birth rate $b$
 and initial condition $\eta _0 (k) = \1 \{k = 0\}$, $k \in \Z$.
  For an interval $[a,b]\subset \R$ and $c > 0$, $c[a,b] $ denotes the interval $[ca, cb].$
 The following theorem characterizes the growth of the set of occupied sites.
 \begin{thm}\label{shape thm}
 	There exists $\lambda _r, \lambda _l >0$ such that for every $\varepsilon >0$ a.s. for sufficiently large $t$,
 	\begin{equation}\label{alcove}
   \bigg(	(1 - \varepsilon)[-\lambda _l t, \lambda _r t] \cap \Z \bigg) \subset \textnormal{occ}(\eta _t)
 	\subset (1 + \varepsilon)[-\lambda _l t, \lambda _r t].
 	\end{equation}
 \end{thm}

\begin{rmk}
	If we assume additionally that the birth rate is symmetric, that is,
    if	for all $x \in \Z$, $\eta \in \X$, 
	\[
	b(x, \eta ) = b(-x,  \tilde \eta),
	\]
	where $\tilde \eta (y) = \eta (-y)$,
	then, as can be seen from the proof, $\lambda _l = \lambda _r$ holds true in Theorem 
	\ref{shape thm}.
\end{rmk}

\begin{rmk}\label{rmk attractiveness}
	Note that
	under our assumptions
	the following attractiveness property does not have to hold: if for two initial configuration
	$\eta ^1 _0 \leq \eta ^2 _0$, then $\eta ^1 _t \leq \eta ^2 _t$ for all $t\geq 0$. This
	renders inapplicable the techniques based 
	on a  subadditive ergodic theorem (e.g. \cite{Lig85subadd})
	which are usually used in the proof of shape theorems (see e.g. \cite{Dur88, Dei03, shapenodeath}).
	On the other hand,
	our technique  relies heavily on
	the dimension being one as
	the  analysis is based on 
	viewing the process from its tip.
	It would be of interest to extend the 
	result to dimensions $\d \geq 2$.
	To the best of our knowledge, even for the following
	modification of Eden's model, 
	the shape theorem has not been proven.
	Take $\d = 2$, $N = 1$ (only one particle per site is allowed),
	 and for $x \in \Z ^2$ and
	$\eta \subset \Z ^2$ with $\eta (x) = 0$ let 
	\[
	b(x, \eta) = \1 \left\{  \sum\limits _{y \in \Z ^2: |y - x| = 1}\eta (y) \in \{1,2\}  \right\}.
	\]
	and define $\xi _t = \eta _t \cup \{y \in \Z ^2: y \text{ is surrounded by particles of } \eta _t  \}$, $t \geq 0$.
	It is reasonable to  expect the shape theorem to hold
	for  $\xi _t$.
	Note that the classical Eden model can be seen as a birth process with rate 
	$$b(x, \eta) = \1\left\{  \sum\limits _{y \in \Z ^2: |y - x| = 1}\eta (y) > 0  \right\}$$
	started from a single particle at the origin.
\end{rmk}
 
For $\eta \in \X$ with $\sum\limits _{x \in \Z} \eta (x) < \infty$, let 
$$\text{tip}(\eta) := \max\{m \in \Z: \eta (m) > 0\}$$
be the position of the rightmost occupied site. Let 
\begin{equation}\label{X = tip}
X_t = \text{tip} (\eta _t).
\end{equation}
By Theorem \ref{shape thm}  a.s. 
\begin{equation}\label{protrude}
 \frac{X_t}{t} \to \lambda _r,  \  \ \ t \to \infty.
\end{equation}

We now give two results on the deviations  of $X_t$ from the mean.
The first theorem gives a sub-Gaussian  limiting estimate on the fluctuations around the mean,
while the second provices an exponential estimate for all $t \geq 0$.
Let $\lambda _r$ be as in  \eqref{protrude}.

\begin{thm}\label{fluctuations thm Gauss est}
	
	There exist $\mathrm{C}_1, \sigma  >0$ such that  
	\begin{equation}
	\limsup\limits_{t \to \infty}
	\PP{ |X_t - \lambda _r t| \geq q \sqrt{t} } \leq \mathrm{C}_1
	 e^{-q ^2 \sigma ^2},
	 \ \ \ q > 0
	.
	\end{equation}

\end{thm}

\begin{thm}\label{fluctuations thm exp est}
     There exist $\mathrm{C}_2, \theta  >0$ such that
     	\begin{equation}\label{putrid}
     \PP{ |\frac{X_t}{t} - \lambda _r | \geq  q } \leq \mathrm{C}_2 e^{-\theta q t},
     \ \ \ q > 0, \ t > 0
     .
     \end{equation}
\end{thm}

Of course, Theorems \ref{fluctuations thm Gauss est} and \ref{fluctuations thm exp est} also apply to 
the position of the leftmost occupied site provided that $\lambda_r$
is replaced with $\lambda _l$.

\label{est fec page}
\emph{Birth rate with regulation via fecundity and establishment}.
As an example of a non-trivial model satisfying our assumptions,
consider the birth process in $\X$ with
birth rate 
\begin{equation}\label{fec est}
\begin{multlined}
b(x, \eta) = \exp \left\{ - \sum\limits _{u \in \Z } \phi(u-x) \eta (u) \right\} 
\sum\limits _{y \in \Z } \left[ a(x-y) \eta (y) 
 \exp \left\{ -\sum\limits _{v \in \Z } \psi(v-y) \eta (v)  \right\} \right], 
 \\
  x\in \Z , \eta \in \X, \eta (x) < N. 
\end{multlined}
\end{equation}
where $a, \phi, \psi: \Z \to \R _+$ 
have a finite range, $\sum\limits _{x \in \Z} a(x) > 0$. The birth rate \eqref{fec est}
is a modification of the free branching rate 
\begin{equation*}
\begin{multlined}
b(x, \eta) =\sum\limits _{y \in \Z }  a(x-y) \eta (y) ,  \ \ \ 
x\in \Z , \eta \in \X, \eta (x) < N. 
\end{multlined}
\end{equation*}

The purpose of the modification is to include damping mechanisms
reducing the birth rate  
in the dense regions. The first exponent multiplier 
in \eqref{fec est}, $ \exp \left\{ - \sum\limits _{u \in \Z } \phi(u-x) \eta (u) \right\}$,
represents the reduction in establishment 
at location $x$ if $\eta$ has many individuals around $x$. 
The second exponent multiplier, $\exp \left\{ -\sum\limits _{v \in \Z } \psi(v-y) \eta (v)  \right\}$,
represents
diminishing fecundity of an individual {at $y$} surrounded by many other individuals.
Further description and
motivation for
 an equivalent continuous-space model 
 can be found in \cite{est_fec, fec_19}.
We note here that the birth process with birth rate \eqref{fec est}
does not in general possess the attractiveness property mentioned in Remark \ref{rmk attractiveness}.
Some numerical observations on this model are collected in Section \ref{sec: numerics}.

\section{Construction of the process}\label{sec: construction}

 Similarly to \cite{Garcia}, \cite{GarciaKurtz}, \cite{shapenodeath}, we construct the process as a solution to 
 the stochastic equation
 
 \begin{equation} \label{se}
 \begin{aligned}
 \eta _t (k) = \int\limits _{(0,t] \times \{k\} \times [0, \infty )   }
  \mathds{1} _{ [0,b(i, \eta _{s-} )] } (u) 
 P(ds,di,du)
 + \eta _0 (k),  \quad
 t \geq q, \  k \in \Z,
 \end{aligned}
 \end{equation} 
where
$(\eta _t)_{t \geq 0}$ is a c\`adl\`ag $\X$-valued solution
process,
$P$ is a  Poisson point process on 
$\R_+  \times \Z \times \R_+  $,
the mean measure of $P$ is $ds \times \# \times du$ ($\#$ is the counting measure
on $\Z$ ).
We require the
processes $P$ and $\eta _0$ to be independent of each other. 
Equation \eqref{se} is understood in the sense that the equality
holds a.s. for every 
$k \in \Z $ and $t \geq 0$.
In the integral on the right-hand side of \eqref{se},
$i = k$ is the location and 
$s$ 
is the time of birth of a new particle. Thus, 
the integral  from $0$ to $t$ represents the number 
of births at $k \in \Z$ which occurred before $t$.

This section follows closely 
Section 5 in \cite{shapenodeath}.
Note that the only difference to Theorem 5.1 from \cite{shapenodeath}
is that the `geographic' space  is discrete ($\Z$) rather than continuous ($\R ^\d$ as in \cite{shapenodeath}).
This change  requires no new arguments,  ideas, or techniques  in comparison to \cite{shapenodeath}.

 We will make the following assumption on the initial  condition:
 \begin{equation} \label{condition on eta _0}
 \E \sum\limits _{i \in \Z}\eta _0 (i) < \infty.
 \end{equation}
 
 Let $P$ be defined 
 on a probability space $(\Omega , \mathscr{F} , \P) $.
 We say that the process
 $P$ is  \textit{compatible} with 
 an increasing, right-continuous 
 and complete filtration of 
 $\sigma$-algebras
 $(\mathscr{F}_t, t \geq 0)$, $\mathscr{F}_t \subset \mathscr{F}$,
 if $P$ is adapted, that is, all random variables
 of the type $P(\bar T _1 \times U)$, $\bar T _1 \in \mathscr{B}([0;t])$,
 $U \in \mathscr{B}(\Z \times \R _+)$, are $\mathscr{F}_t$-measurable,
 and all random variables of the type
 $P( (t , t + h] \times U) $, $ h \geq 0$, $U \in \mathscr{B}(\Z \times \R _+)$,
 are independent of $\mathscr{F}_t$ (here we consider 
 the Borel $\sigma$-algebra for $\Z$ to be the collection  $2 ^{\Z}$ of all subsets of $\Z$).

We equip $\X$ with the product set topology and $\sigma$-algebra generated by the open sets in this topology.
 
 \begin{defi} \label{weak solution}
 	A \emph{(weak) solution} of equation \eqref{se} is a triple 
 	$(( \eta _t )_{t\geq 0} , P )$, $(\Omega , \mathscr{F} , \P) $,
 	$(\{ \mathscr {F} _t  \} _ {t\geq 0}) $, where

 (i) $(\Omega , \mathscr{F} , \P)$ is a probability space, 
 and $\{ \mathscr {F} _t  \} _ {t\geq 0}$ is an increasing, right-continuous
 and complete filtration of sub-$\sigma$-algebras of $\mathscr {F}$,
 
 (ii)   $P$ is a
 Poisson point process on $\R_+  \times \Z \times \R_+ $  with  
 intensity  $ds \times \# \times du $,

 (iii) $ \eta _0  $ is a random  $\mathscr {F} _0$-measurable
 element in $\X$
 satisfying \eqref{condition on eta _0},
 
 (iv) the processes $P$ and $\eta _0$ are independent,
 $P$ is
 compatible with $\{ \mathscr {F} _t  \} _ {t\geq 0} $,
 
 (v) $( \eta _t )_{t\geq 0} $ is a c\`adl\`ag $\X$-valued process
 adapted to $\{ \mathscr {F} _t  \} _ {t\geq 0} $, $\eta _t \big| _{t=0} = \eta _0$,  
 
 (vi) all integrals in \eqref{se} are well-defined, 
 \[
 \E \int\limits _0 ^t ds \sum\limits _{i \in \Z} b(i, \eta _{s-}) 
 < \infty, \ \ \ t > 0,
 \]
 
 (vii) equality \eqref{se} holds a.s. for all $t\in [0,\infty]$
 and all $k \in \Z$.
\end{defi}

 Let 
 \begin{align}\label{nonchalant}
 \mathscr{S} ^{0} _t =  \sigma \bigl\{ & \eta_0 , 
 P([0,q] \times \{k\} \times C ) ,   \\ &
 q \in [0,t], 
 k \in \Z, C \in \mathscr{B} (\R_+)  \bigr\},
 \notag
 \end{align}
 and
 let $\mathscr{S} _t$ be the completion of $\mathscr{S} ^{0} _t$ under $P$.
 Note that $\{ \mathscr{S} _t \}_{t\geq 0} $ 
 is a right-continuous filtration.

 \begin{defi} \label{def strong solution}
 	A solution  of \eqref{se} is called \emph{strong}
 	if $( \eta _t )_{t\geq 0} $ is adapted to 
 	$(\mathscr{S} _t, t\geq 0)$.
 \end{defi}

 \begin{defi}
 	We say that \emph{pathwise uniqueness} holds for \eqref{se} if for any two (weak)
 	solutions 	$(( \eta _t )_{t\geq 0} , P )$, $(\Omega , \mathscr{F} , \P) $,
 	$(\{ \mathscr {F} _t  \} _ {t\geq 0}) $ 
 	and 
 	$(( \eta _t ^{\prime} )_{t\geq 0} , P )$, $(\Omega , \mathscr{F} , \P) $,
 	$(\{ \mathscr {F} _t  \} _ {t\geq 0}) $ 
 	with $\eta _0 = \eta _0 ^{\prime}$ we have 
 	\begin{equation}
 	\PP{\eta _t = \eta _t ^{\prime} \text{ for all } t\geq 0} = 1.
 	\end{equation}
 \end{defi}

\begin{defi} \label{joint uniqueness in law}
	We say that \textit{joint uniqueness in law} holds for equation \eqref{se} with an initial
	distribution $\nu$ if any two (weak) solutions $((\eta_t) , P)$ and 
	$((\eta_t  ^{ \prime }) , P  ^{\prime} )$ of \eqref{se},
	$Law(\eta _0)= Law(  \eta _0  ^{\prime})=\nu$, have the same joint distribution:
	
	$$Law ((\eta_t) , P)
	= Law ((\eta_t  ^{\prime} ), P ^{\prime}  ).$$

\end{defi}
 
 \begin{thm}\label{core thm}
 	Pathwise uniqueness, strong existence and joint uniqueness
 	in law hold for equation \eqref{se}. 
 	The unique solution is a  Markov process with respect to 
 	the filtration $(\mathscr{S} _t, t\geq 0)$. 
 \end{thm}
 
 The proof follows exactly the proof of Theorem 5.1 in \cite{shapenodeath} and  is therefore omitted.
 We also note here that  the unique solution of \eqref{se} satisfies a.s. $\eta _t(x) \leq N$
 for $x \in \Z$ and $t \geq 0$ by \eqref{winnow}. 
 
 \section{Proofs} \label{sec: proof}

 Let
 $\Z _+ = \N \cup \{0\}$, $\Z _- = -\Z _+$, and
  $\beta _t: \Z _- \to \Z _+$ be $\eta _t$ seen from its tip
 defined by
 \[
  \beta _t (-n) = \eta _t (\text{tip}(\eta _t) - n), \ \ \ n = 0,1,2,...
 \]
 Let $h_t$ be the position of first block of $R$ sites occupied by $N$ particles
 seen from the tip,
 \[
 h_t := \max\{m \in \Z _-: \beta _t(m-1) = \beta _t(m-2) = ... = \beta _t(m-R) = N   \}
 \vee \min\{ m \in \Z _-: \beta _t (m) >0 \}.
 \]
 We adopt here the convention $\max\{ \varnothing\} = -\infty$, so that 
 if there are no blocks of  $R$ consecutive sites occupied by $N$ particles,
 $h_t$ equals to the furthest from the origin occupied site for $\beta _t$.
 Finally, define $\alpha  _t: \Z _- \to \{0,1,...,N\} $ by
 \[
 \alpha _t (m)= \beta _t (m) \I \{m \geq h_t\}.
 \]
 Thus, $\alpha _t$ can be interpreted as the part of
  $\eta _t$ seen from its tip until the first block of $R$ sites occupied by $N$
  particles. 
  The process $(\alpha _t, t \geq 0)$ takes values in 
  a countable space 
  \[
  \Upsilon :=  \left\{ \gamma \middle| \gamma : \Z _- \to \{0,1,...,N\}, \sum\limits_{x \in \Z_-}\gamma (x) < \infty  \right\}.
  \]
  Let us underline that  $\alpha _t$ is a function of  $\eta _t$;
  we denote by $\mathcal{A}$ the respective mapping 
  $\mathcal{A} : \X \to \Upsilon $,
  so that
   $\alpha _t = \mathcal{A}(\eta _t)$.
\label{mathcal A}

  \begin{lem}\label{gyrate}
  	The process
  	$(\alpha _t, t \geq 0)$ is a continuous-time positive recurrent Markov process
  	with a countable state space. Furthermore, $(\alpha _t, t \geq 0)$ is strongly ergodic.
  \end{lem}
 
 	\textbf{Proof}. 
 	We start from 
 	 a key observation: for $t\geq 0$,
 	conditionally on the event  
 	$$\eta (m) = \eta (m+1) = ... = \eta (m+R) = N$$
 	for some $m \in \Z$,
 	the families $\{\eta _s(y), s\geq t, y < m  \}$ and $\{\eta _s (y), s\geq t, y > m +R \}$
 	are independent. 
    Consequently, $(\alpha _t, t \geq 0)$ is an irreducible
    continuous-time Markov chain.
    The definitions and properties of continuous-time Markov chains
    used here
    can be found e.g. in \cite[Section 4.4]{ChenMarkov}.
    Translation invariance Condition \ref{con translation invariance}
    ensures that $(\alpha _t, t \geq 0)$ is time-homogeneous.
      Define $\0 _{\Upsilon} \in \Upsilon$  by $\0 _{\Upsilon} (m) = 0 $, $m = 0,-1,-2,...$
	Now, let us note that
	
	\begin{equation}\label{morass}
	 \inf\limits _{\gamma \in \Upsilon} \Pc{\alpha _{t+1} = \0 _{\Upsilon}}{\alpha _t = \gamma} 
	 =
	 \inf\limits _{\gamma \in \Upsilon} \Pc{\alpha _{1} = \0 _{\Upsilon}}{\alpha _0 = \gamma}  >0.
	\end{equation}
	Indeed, for $\eta _0$ such that $\mathcal{A}(\eta_0) = \alpha _0 = \gamma$,
	\[
	\PP{\alpha _{1} = \0 _{\Upsilon}} \geq \PP{ \text{tip} (\eta _1) = \text{tip} (\eta _0)}
	\]
	\[
	\times \PP{\eta_1(\text{tip} (\eta _1)) = \eta_1(\text{tip} (\eta _1-1)) = ... = \eta_1(\text{tip} (\eta _1)-R+1) = N} >0,
	\]
	and the last expression is separated from $0$ uniformly in $\gamma$.
	It follows from \eqref{morass} that the state $\0 _{\Upsilon}$
	is positive recurrent. Since $(\alpha _t, t \geq 0)$ is irreducible,
	it follows that it is also positive recurrent.
	The strong ergodicity follows from \eqref{morass}.
	\qed

	Denote by $\pi $ the ergodic measure for $(\alpha _t, t \geq 0)$.
	For $\gamma \in \Upsilon$, let $\eta ^\gamma \in \X$ be
	\begin{equation}
	\eta ^\gamma (m) = \begin{cases}
						 \gamma (m), &  \ \ \ m \leq 0 
						 \\
						 0, & \ \ \ m >0
				       \end{cases}
	\end{equation}
	Note that $\mathcal{A}(\eta ^\gamma ) = \gamma$.
	Define $f:\Upsilon \to \R _+$ by
	\begin{equation*}
	f(\gamma ) = \sum\limits _{m = 1} ^N m b(x, \eta ^\gamma).
	\end{equation*}
	
	Note that 
	\begin{equation}
	\sup\limits _{\gamma \in \Upsilon} f(\gamma) \leq \frac{R(R + 1)}{2}  \overline{\mathbf{ b}} .
	\end{equation}
	Since $f$ is bounded,
	 by the ergodic theorem for continuous-time 
	Markov chains
	a.s.
	
	\begin{equation}\label{condiment}
	\frac 1t \int\limits _{0} ^t f(\alpha _s)ds \to \langle f \rangle _{\pi},
	\end{equation}
 	where $\langle f \rangle _{\pi} := \sum\limits _{\gamma \in \Upsilon} \pi (\gamma) f(\gamma)$ (here for convenience 
 	$\pi(\{\gamma\})$ is denoted by $\pi(\gamma)$, $\gamma \in \Upsilon$).

	Recall that $X_t = \text{tip}(\eta _t)$. The process  $(X_t, t \geq 0)$ is an increasing pure jump type Markov process,
	and  the rate of jump of size $m \in \{1,2,...,R\}$
	at time $t$ is  $b(X_t + m, \eta_t)$.
	Indeed,
	note that 
	\begin{equation}\label{visceral}
	\begin{split}
	X_t & = \sum\limits_{k = 1} ^R  k
	\int\limits _{(0,t] \times \Z \times [0, \infty )   } 
	\mathds{1} \{ i = k +X_{s-}  \}
	\mathds{1} _{ [0,b(k + X _{s-}, \eta _{s-} )] } (u) 
	P (ds,di,du) 
	\\
    &	=  \sum\limits_{k = 1} ^R  k
		\int\limits _{(0,t]  \times [0, \infty )   } 
	\mathds{1} _{ [0,b(k + X _{s-}, \eta _{s-} )] } (u) 
	P^{(k)} (ds,du),
	\end{split}
	\end{equation}
	where the integrator is defined by 
	\begin{equation*}
	\begin{multlined}
	P^{(k)} (A \times B) = P\big (\{(t, i, u) \in \R _+ \times \Z \times \R _+ \big| X_{t-} + k = i, (t,u) \in A \times B \} \big),
	\\
	\
	\ A, B \in \mathscr{B}( \R _+), \ \ k \in \{1,...,m\}.
		\end{multlined}
	\end{equation*}
   In other words, $P^{(k)}(A \times B) $ is $ P(A \times \{i\} \times B)$ if  $ A \subset \{t: X_{t-} + k = i\}$
   and $B \in \mathscr{B}(\R _+)$.
   Note that $P^{(k)}$ 
	is a Poisson point process on $\R _+ \times \R _+$ with mean measure $ds \times du$ (this follows
	for example from the strong Markov property of a Poisson point process, as formulated in the appendix
	in \cite{shapenodeath}, applied to the jump times of $(X_t, t\geq 0)$ ). 
	The indicators in \eqref{visceral} are
	\begin{alignat*}{2}
	& \begin{aligned} & 
	\mathds{1} \{ i = k +X_{s-}  \}
	=
	\begin{cases}
	1, & \text{ if }  i = k +X_{s-}, \\
	0, & \text{ otherwise },\\
	\end{cases}\\
	\end{aligned}
	& \hskip 2em &
	\begin{aligned}
	&
	\mathds{1} _{ [0,b(k + X _t, \eta _{s-} )] } (u) =
	\begin{cases}
    1, & \text{ if }  u \in [0,b(k + X _t, \eta _{s-} )], \\
     0, & \text{ otherwise }.\\
     \end{cases}\\
	\end{aligned}
	\end{alignat*}

	Therefore, by e.g. (3.8) in Section 3, Chapter 2 of \cite{IkedaWat}, the process 
	\begin{equation}\label{martingale}
	   M _t : = X_t - \int\limits _0 ^t \sum\limits _{k=1} ^R k b(X_t + k, \eta_s) ds = 
	   X_t - \int\limits _0 ^t f(\alpha _s)ds.
		\end{equation}
	is a martingale with respect to the filtration $(\mathscr{S} _t, t\geq 0)$ defined below \eqref{nonchalant}.
	
	We now formulate a strong law of large numbers for martingales.
        The following theorem is an abridged version of \cite[Theorem 2.18]{HH80}.

\begin{thm}\label{thm 2.18 HH80}
   Let $\{S_n = \sum\limits _{i = 1} ^n x _i , n \in \N  \}$ be an $\{\mathscr{S}_n\}$-martingale
   and $\{U_n\}_{n \in \N }$ be a non-decreasing sequence of positive real numbers, $\lim\limits_{n \to \infty} U_n = \infty$.
   Then for $p \in [1,2]$ we have
   \[
   \lim\limits_{n \to \infty } U _n ^{-1} S _n = 0
   \]
   a.s. on the set $\left\{ \sum\limits _{i = 1} ^\infty U _n ^{-p} \E \big\{|x_i|^p \big | \mathscr{F}_{i-1}\big\} < \infty  \right\}$.
\end{thm}
	
	\begin{lem}\label{put on airs}
		Strong law of large numbers  applies to $(M_t, t \geq 0)$: 
		\begin{equation}\label{denigrate}
		 \PP{\frac{M_t}{t} \to 0} = 1.
		\end{equation}
	\end{lem}
	\textbf{Proof}. Let $\Delta M _n = M_{n+1} - M_n$.
	Then  $\Delta M_n$
	is stochastically 
	dominated by $\mathrm{V}_1 + 2 \mathrm{V}_2 + ... + R\mathrm{V}_R + \frac{R(R + 1)}{2}  \overline{\mathbf{ b}} $,
	where
	$\mathrm{V}_1,..., \mathrm{V}_R$ are independent Poisson random variables
	with mean $\overline{\mathbf{ b}}$, independent of $\mathscr{F}_n$. 
	Applying the strong law of large numbers for martingales from
	Theorem \ref{thm 2.18 HH80} with $p = \frac 32$ and $U_n = n$, we get a.s.
	\begin{equation}
	\frac 1n \sum\limits \Delta M_n \to 0.
	\end{equation}
	Since a.s. for every $\varepsilon >0$, 
	$$
	\PP{\sup\limits _{s \in [0,1] }|M_{n+s} - M_s| \geq \varepsilon n
	\text{ for infinitely many } n}
    = 0, $$
	 \eqref{denigrate} follows.
	\qed
	
	\textbf{Proof of Theorem \ref{shape thm}}.
	Let $\lambda _r = \langle f \rangle _{\pi}$.
	From \eqref{condiment} and Lemma \ref{put on airs} we get  a.s.
	\begin{equation}
	\frac 1t  X_t - \langle f \rangle _{\pi} \to 0,
	\end{equation}
	or 
		\begin{equation}\label{flip out}
	\frac 1t  \text{tip}(\eta _t) -\lambda _r \to 0.
	\end{equation}
	In the same way (due to the symmetric nature of our assumptions)
	we can show the equivalent of \eqref{flip out} for the leftmost occupied
	site $Y_t$ for $\eta _t$: there exists $\lambda _l >0$
	such that
		for any $\varepsilon >0$ a.s. for large $t$, 
		\begin{equation}\label{flip out left}
		\frac { |Y_t|}{t} - \lambda _l \to 0.
		\end{equation}
	Hence the second inclusion in \eqref{alcove} holds.
 (As an aside we point out here that $ \lambda _l$ 
 can be expressed as an average value in the same way as
  $\lambda _r = \langle f \rangle _{\pi}$. To do so, we would need
  to define the opposite direction counterparts   
  to $(\beta _t, t \geq 0)$, 
  $(\alpha _t, t \geq 0)$, $f$, and other related objects.)

	To show the first inclusion in \eqref{alcove},
	we fix $\varepsilon >0$, and for each $x \in \Z$. 
	By \eqref{flip out} and \eqref{flip out left},
	a.s. for large $t$ 
	\begin{equation}\label{connoisseur}
	 \left[-\lambda _l t \left(1 + \frac \varepsilon 4 \right)^{-1}, \lambda _r t \left(1 + \frac \varepsilon 4\right)^{-1}\right] \subset [Y_t, X_t].
	\end{equation}	
	For $x \in \Z$, let 
	\begin{equation}
	\begin{split}
	\sigma (x) & := \inf\{t\geq 0: dist(x, \text{occ}(\eta _t)) \leq R \},
		\\
	\tau (x) & := \inf\{t\geq 0: \eta_t (x) \geq 1 \} = \inf\{t\geq 0: x \in \text{occ}(\eta _t)\}.
	\end{split}
	\end{equation}
	Clearly, for any $x \in \Z$, $0 \leq \sigma (x) \leq \tau (x)$. 
	Because of the finite range assumption,
	by \eqref{connoisseur} a.s. for $x \in \N$ with large  $|x|$
	\begin{equation} \label{slink}
	\sigma (x)  \leq \frac{(1 + \frac \varepsilon 4) |x|}{\lambda _r}.
	\end{equation}
	By \eqref{grovel}, the random variable $\tau (x) - \sigma (x)$
	is stochastically dominated by an exponential random variable
	with mean $\underline{\mathbf{b} } ^{-1}$. In particular
	\begin{equation*}
	\PP{ \tau(x) - \sigma(x) \geq \frac{\varepsilon |x|}{4 \lambda _r}} \leq \exp\{ - \frac{\varepsilon |x|\underline{\mathbf{b} }}{4\lambda _r} \}.
	\end{equation*}
	Since $\sum\limits _{x \in \N} \exp\{ - \frac{\varepsilon |x|\underline{\mathbf{b} }}{4\lambda _r} \} < \infty$,
	a.s.
	for all but finitely many $x \in \N$ we have 
	\begin{equation*}
	\tau(x) - \sigma(x) \leq \frac{\varepsilon |x|}{4\lambda _r}.
	\end{equation*}
	Hence from \eqref{slink} a.s. for all but finitely many $x \in \N$,
	\begin{equation}\label{bland}
	 \tau (x ) \leq \frac{(1 + \frac \varepsilon 2) |x|}{\lambda _r}.
	\end{equation}
	From \eqref{bland} it follows that a.s. if $t$ is large 
	and $x \in \N$, $|x| \leq\frac{\lambda _r t}{1 + \frac \varepsilon 2} $, 
	then $\tau (x) \leq t$. Note that $\text{occ}(\eta _t ) = \{ x\in \Z: \tau (x) \leq t\}$. Thus
	for large $t$. 
	\begin{equation}\label{faucet}
    \left[0,  \frac{\lambda _r t}{1 + \frac \varepsilon 2}\right] \cap \Z  \subset 	\text{occ}(\eta _t).
	\end{equation}
	
	Repeating this argument
	verbatim
	 for $-x \in \N$ and $\lambda _l$ in place of $x \in \N$
	and $\lambda _r$ respectively, we find 
	that 
		\begin{equation}\label{quail}
	\left[ - \frac{\lambda _l t}{1 + \frac \varepsilon 2}, 0\right]  \cap \Z  \subset 	\text{occ}(\eta _t).
	\end{equation}
	Since $ 1 - \varepsilon \leq \frac{1}{1 + \frac \varepsilon 2}$ for $\varepsilon >0$,
	the first inclusion in  \eqref{alcove} 
	follows from \eqref{faucet} and \eqref{quail}.
	\qed

	\begin{lem}\label{sulky}
		For some $\varrho ^2 \in (0, +\infty)$
		a.s.
		\begin{equation}\label{deluge}
		\frac 1t  [M]_t \to \varrho ^2, \ \ \ t\to \infty.
		\end{equation}
	\end{lem}
	\textbf{Proof}.
	Let $\theta _0 = 0$  and
	denote by $\theta _n$,   $n \in \N$,
	the moment of $n$-th hitting of 
	$\0 _\Upsilon$
	by the Markov chain
		$(\alpha _t, t \geq 0)$.
		For $n \in \N$, define
		a random piecewise constant  function       
		 $Z_n $
		 by 
		 \begin{equation}
		 Z_n (t)= \alpha _{  (t + \theta _n) \wedge \theta _{n+1} },  \ \  \ t\geq 0.
		 \end{equation}
	The sequence $\{Z_n\}_{n \in \N}$ can be seen as a sequence of independent 
	random elements in  the Skorokhod space $\mathcal{D}: = D([0,\infty), \Upsilon)$
	endowed with the usual Skorokhod topology.
	
	Let $G: \mathcal{D} \to \Z _+$ be the functional
	such that $G(Z_n) =  [M]_{\theta _{n+1} - \theta _{n}}$
	  is the change of $( [M] _t)$
	 between $\theta _n$ and $\theta _{n+1}$. The function $G$ can be written down explicitly,
	 but it is not necessary for our purposes.
	Now, since the number of jumps for $Z_1$ has exponential tails, 
	for any $m \in \N$
	\begin{equation}
	  \EE{G^m(Z_1)} < \infty.
	\end{equation}
	
	By the strong law of large numbers, a.s.
	\begin{equation}\label{0909-1}
	\frac 1 n \sum\limits _{i = 1} ^n G(Z_i) \to \EE{G(Z_1)} >0.
	\end{equation}
	
	Since $\theta _2 - \theta _1$, $\theta _3 - \theta _2$, ..., are i.i.d. random variables, a.s.
	\begin{equation}\label{0909-2}
	 \frac{\theta _n}{n} \to \EE{\theta_2 - \theta _1} >0. 
	\end{equation}

	For $t>0$, let $n = n(t) \in \N$  be such that $t \in [\theta _n, \theta _{n+1})$. 
	Then  by  \eqref{0909-1} and \eqref{0909-2} a.s.
	\begin{equation}
	\lim\limits _{t \to \infty} \frac{[M]_t}{t}
	=
	\lim\limits _{t \to \infty}\frac{[M]_{\theta _1} +  \sum\limits _{i = 1} ^{n(t)} G(Z_i)}{t}
	=
	\lim\limits _{t \to \infty}\frac{\sum\limits _{i = 1} ^{n(t)} G(Z_i)}{n(t)}
	\frac{n(t)}{t}
	= \frac{\EE{G(Z_1)}}{\EE{\theta_2 - \theta _1}} >0.
	\end{equation}
\qed

	Before proceeding with the final part of the paper, 
	we formulate a central limit theorem for martingales used
	in
	the proof of Theorem \ref{fluctuations thm Gauss est}.
	The statement below is a corollary  of \cite[Theorem 5.1]{CLT_mart_CT}.
	
	\begin{prop}\label{thm 5.1 CLT_mart_CT}
	 Assume that \eqref{deluge} holds and that 
	 for some $K >0$ a.s.
	 $$
	 \sup\limits _{t \geq 0} |M_t - M_{t-}| \leq K.
	 $$
	 Then 
	 	\begin{equation}
	\frac{1}{\sqrt{t}}   M _t
	\overset{d}{ \to } 
	\mathcal{N}(0, \varrho ^2 ), \ \ \ t\to \infty.
	\end{equation}
	\end{prop}
      Proposition \ref{thm 5.1 CLT_mart_CT} follows from \cite[Theorem 5.1, (b)]{CLT_mart_CT}
      by taking  $M_n$ in notation 
      of \cite{CLT_mart_CT}
      to be $\frac{M_{tn}}{n}$ in our notation.
	
	\textbf{Proof of Theorem \ref{fluctuations thm Gauss est}}.
  By Lemma \ref{gyrate},	
	the continuous-time Markov chain $(\alpha _t, t \geq 0)$
	is strongly ergodic.
	  Since the function $f$ is bounded,
	the central limit theorem holds
	for $(\alpha _t, t \geq 0)$
		by \cite[Theorem 3.1]{CLT_for_CTMC}.
		That is, the convergence in distribution takes place
		\begin{equation}\label{discombobulate}
		\frac{1}{\sqrt{t}}  \left[\int\limits _0 ^t f(\alpha _s )ds  - \langle f \rangle _{\pi} 
		 \right]
		 \overset{d}{ \to } 
		 \mathcal{N}(0, \sigma ^2 _{f}), \ \ \ t\to \infty,
		\end{equation}
	where $\sigma ^2 _{f} \geq 0$ is a constant depending on $f$,
	and $\mathcal{N}(c, \sigma ^2 )$ is 
	the normal distribution with mean $c\in \R$ 
	and  variance $\sigma ^2 \geq 0$. Recall that 
	$\langle f \rangle _{\pi}  = \lambda _r$.
	By \eqref{discombobulate},
	for some $C_1, \sigma _1 >0$
		\begin{equation}\label{subgaus1}
	\limsup\limits_{t \to \infty}
	\PP{ \left| \int\limits _0 ^t f(\alpha _s )ds - \lambda _r t \right| \geq q \sqrt{t} } \leq C_1 e^{-q ^2 \sigma _1 ^2},
	\ \ \ q > 0
	.
	\end{equation}
	
	Recall that $M_t$ was defined in \eqref{martingale}.
	By Lemma \ref{sulky}  for some $\varrho ^2 \in (0, +\infty)$
	a.s.
	\begin{equation}
	\frac 1t  [M]_t \to \varrho ^2, \ \ \ t\to \infty.
	\end{equation}
	By the martingale central limit theorem (Proposition \ref{thm 5.1 CLT_mart_CT})
		\begin{equation}
	\frac{1}{\sqrt{t}}   M _t
	\overset{d}{ \to } 
	\mathcal{N}(0, \varrho ^2 ), \ \ \ t\to \infty.
	\end{equation}
	Hence 
		for some $C_2, \sigma _2 >0$
	\begin{equation} \label{subgaus2}
		\limsup\limits_{t \to \infty}
	\PP{ \left| X_t -  \int\limits _0 ^t f(\alpha _s )ds \right| \geq q \sqrt{t} } \leq C_2 e^{-q ^2 \sigma _2 ^2},
	\ \ \ q > 0
	.	
	\end{equation}
	By \eqref{subgaus1} and \eqref{subgaus2},
		\begin{equation} \label{subgaus3}
	\limsup\limits_{t \to \infty}
	\PP{ \left| X_t -  \lambda _r t \right| \geq q \sqrt{t} } \leq C_3 e^{-q ^2 \sigma _3 ^2},
	\ \ \ q > 0,
	\end{equation}
	for some $C_3, \sigma _3 >0$.
	\qed

  \textbf{Proof of Theorem \ref{fluctuations thm exp est}}.
  By \eqref{visceral}, \eqref{martingale},
  and \cite[(3.9), Page 62, Section 3, Chapter 2]{IkedaWat},
  the predictable quadratic variation
  \begin{equation}\label{pred quad var}
   \langle M \rangle _t = \langle X \rangle _t = 
   \int\limits _{0} ^t \sum\limits _{k=1} ^R k ^2 b(k + X_{s-} , \eta _{s-} ) ds
   =
   \int\limits _{0} ^t g(\alpha _{s-})ds,
  \end{equation}  
	where $g: \Upsilon \to \R _+ $ is such that
	 $g(\alpha) = \sum\limits _{k=1} ^R k ^2 b(k + \text{tip}(\eta) , \eta  )$
	whenever $\mathcal{A}(\eta) = \alpha$. 
	Recall that the mapping $\mathcal{A}$
	was defined on Page \pageref{mathcal A}; $g(\alpha)$ does not depend on the choice of $\eta \in \mathcal{A}^{-1}(\alpha)$.
	
	By \cite[Theorem 1.1]{CG08}
	(see also \cite[Theorem 1 and Remark 3a]{Wu00})
   \begin{equation}\label{protruding}
	\P\left\{ \left|\frac 1t \int\limits _{0} ^t g(\alpha _{s-}) ds - \langle g \rangle _{\pi} 
	\right| \geq q \right\} \leq C_1 e ^{-\delta _q t}, \ \ \ q >0, t\geq 0,
    \end{equation}
	where $C_1 >0$,	 $ \delta _q >0$ depends on $q$ but not on $t$,
	and $\delta _q$ 
	grows not slower than linearly as a function of $q$.
	Note  that the jumps of $(M_t, t \geq 0)$
	do not exceed $R$.
	By an exponential inequality for martingales with bounded jumps,
	 \cite[Lemma 2.1]{vdG95},
	 for any $a,b >0$
	
	\begin{equation}
	\PP{|M_t| \geq a, \langle M \rangle _t \leq b \text{ for some } t \geq 0 }
	\leq
	\exp\left\{ -\frac{a^2}{2(aR + b)}  \right\}.
	\end{equation}
	Taking here $a = rt$, $b = \langle g \rangle _\pi (1 + \varepsilon)t $
	for $r>0$, $ \varepsilon \in (0,1)$, we get 
	for $t \geq 0$
	\begin{equation}\label{tack on1}
	 \PP{|M_t| \geq rt, \langle M \rangle _t \leq \langle g \rangle _\pi (1 + \varepsilon)t }
	 \leq 
	  \exp\left\{ -\frac{r^2}{2(rR +  \langle g \rangle _\pi (1 + \varepsilon) )} t \right\}.
	\end{equation}	
	By \eqref{pred quad var} and \eqref{protruding},
	\begin{equation}\label{supple1}
	\PP{ \langle M \rangle _t > \langle g \rangle _\pi (1 + \varepsilon)t  }
	\leq 
	\PP{ \frac 1t |\langle M \rangle _t -  \langle g \rangle _\pi| \geq \varepsilon \langle g \rangle _\pi }
	\leq C_1 e ^{-\delta _{\varepsilon \langle g \rangle _\pi} t}.
	\end{equation}
	By \eqref{tack on1} and \eqref{supple1}
	we get 
	\begin{equation}\label{dry run}
	\begin{gathered}
	\PP{ |M_t| \geq rt} \leq  \PP{|M_t| \geq rt, \langle M \rangle _t \leq \langle g \rangle _\pi (1 + \varepsilon)t }
	+
	\PP{ \langle M \rangle _t > \langle g \rangle _\pi (1 + \varepsilon)t  } 
	\\
	\leq
	\exp\left\{ -\frac{r^2}{2(rR +  \langle g \rangle _\pi (1 + \varepsilon) )} t \right\}
	+
	C e ^{-\delta _{\varepsilon \langle g \rangle _\pi} t}
	\leq C e^{-\delta _1 r t}, \ \ \ t \geq 0,
	\end{gathered}
	\end{equation}
	where $\delta_1 >0$
	does not
	 depend on $r$
	or $t$.

	Recalling the definition of $M$ in \eqref{martingale},
	we  rewrite \eqref{dry run} as
	\begin{equation}\label{travesty}
	\PP{ \left| X_t -  \int\limits _0 ^t f(\alpha _s )ds \right| \geq r } \leq C_2 e ^{-\delta _2 r t} , 
	\ \ \ t \geq 0,
	\end{equation}
	where $C_2, \delta _2 >0$.
		By \cite[Theorem 1.1]{CG08}
	(or \cite[Theorem 1, Remark 3a]{Wu00}),
   \begin{equation}\label{protruding f}
\P\left\{ \left|\frac 1t \int\limits _{0} ^t f(\alpha _{s-}) ds - \langle f \rangle _{\pi} 
\right| \geq r \right\} \leq C_3 e ^{-\delta _3 r t}, \ \ \  t\geq 0,
\end{equation}
where the constant $\delta _3$
does not depend on $r$.
Combining \eqref{travesty} and \eqref{protruding f}
yields   \eqref{putrid} and completes the proof.
	\qed

\begin{rmk}
	We see from the proofs that the non-degeneracy condition
	can be weakened. In particular, 
	\eqref{grovel} can be removed
	if we instead require 
	$(\alpha _t, t \geq 0)$ 
	to be strongly ergodic
	with the ergodic measure satisfying 
	$\langle f \rangle _{\pi} > 0$.
	Of course, $\text{occ} (\eta)$
	in Theorem \ref{shape thm}
	would need to be replaced with the set of sites
	surrounded by $\text{occ} (\eta)$.	
	Some changes in the proof
	would have to be made.
	In particular,  if $\pi (\0 _{\Upsilon}) = 0$, the moments of time $\theta _n$
	in the proof of Lemma \ref{sulky}
	would need to be redefined 
	as the hitting moments 
	of a state $\gamma \in \Upsilon$
	with $\pi (\gamma) > 0$.
\end{rmk}

\begin{rmk}
	 It would be of interest to see 
	 if
	  the finite range condition can be weakened
	  to include interactions decaying exponentially or polynomially fast 
	  with the distance. 
	  If the interaction range is infinite, there is no reason 
	  for   $(\alpha _t, t \geq 0)$ to be a  recurrent  Markov chain,
	  let alone a strongly ergodic one. 
	  It may  be the case however that
	  the process seen from the tip
	   $(\beta _t, t \geq 0)$
	   turns out to possess some kind of  a mixing property,
	   which would enable application of limit theorems.

\end{rmk}

\section{Numerical simulations and monotonicity of the speed} \label{sec: numerics}

We start this section with the following conjecture
claiming that the speed is a monotone functional of the birth rate. 
	Consider
two
birth processes 
$(\eta ^{(1)}_t, t\geq 0)$ and  $(\eta ^{(2)} _t, t\geq 0)$
with
different birth rates $b_1$ and $b_2$, respectively, 
satisfying conditions of Theorem \ref{shape thm}.
Denote by $ \lambda _r ^{(j)}$ the speed 
at which $(\eta ^{(j)}_t, t\geq 0)$ is spreading to the right
in the sense of  Theorem \ref{shape thm}, $j = 1,2$.

\begin{frage}\label{wrong conjecture}

	Assume that for all $x \in \Z $ and $\eta \in \X$ 
	\begin{equation}
	b_1(x, \eta) \leq  b_2(x, \eta).
	\end{equation} 
	Is it always true that   $\lambda _r ^{(1)} \leq \lambda _r ^{(2)}$?
\end{frage}

The answer to
Question \ref{wrong conjecture} is positive
if $b_2$ is additionally assumed to be monotone in the second argument,
that is, if   
$$ b_2(x, \eta) \leq b_2(x, \zeta),
 \ \ \  x \in \Z$$
whenever $\eta \leq \zeta$. Indeed, in this case 
the two birth processes $(\eta ^{(1)}_t, t\geq 0)$ and  $(\eta ^{(2)} _t, t\geq 0)$
with rates $b_1$ and $b_2$ are coupled in such a way that 
a.s. 
$$\eta ^{(1)}_t \leq \eta ^{(2)}_t, 
\ \ \ t \geq 0
$$ (see Lemma 5.1  in \cite{shapenodeath}).
One might think  that the answer is positive in a general case, too. 

It turns out that 
the birth rate with fecundity and establishment regulation
discussed on Page \pageref{est fec page}
links up naturally with Question \ref{wrong conjecture}. 
Let the birth rate $b$ be as in \eqref{fec est} with $R = N = 3$, $a(x) = \1 \{ |x| \leq 3 \}$,
$$\psi (x) = c_{\text{fec}} \left[ \1 \{ |x| =  0 \} + \frac 12 \1 \{ |x| =  1 \} \right] ,
\ \ \ 
\psi (x) = c_{\text{est}} \left[ \1 \{ |x| =  0 \} + \frac 12 \1 \{ |x| =  1 \} \right]. $$

		\begin{figure}[H]
	\centering
	\includegraphics[scale=0.76]{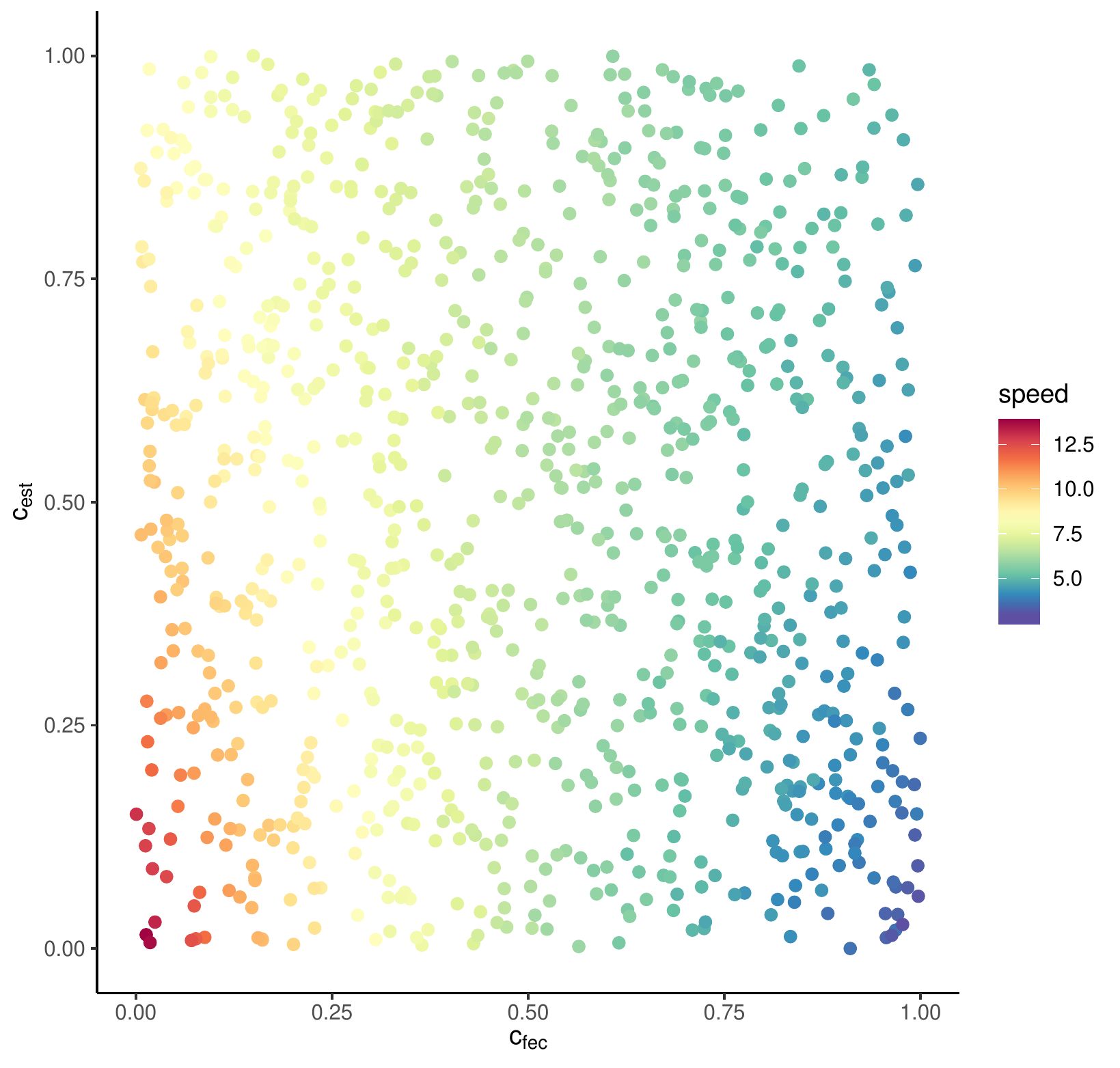}
	\caption{The speed of the tip for various $c_{\text{fec}} , c_{\text{est}}$. 
	For each pair, the speed is computed as the average speed of the tip $X_t$ between
	$t_1 = 100$ and $t_2 = 1000$, that is, as $\frac{X_{t_2} - X_{t_1}}{t_2 - t_1}$.
	Early evolution is excluded to reduce bias.}
	\label{f1}
\end{figure}

Note that $b$ decreases
as either of the parameters  $c_{\text{fec}}$ and $c_{\text{est}}$ increases.
Figure \ref{f1} shows the speed of the model with birth rate \eqref{fec est}
for a thousand randomly chosen from $[0,1]^2$ pairs of parameters $(c_{\text{fec}} , c_{\text{est}})$.

\begin{figure}[H]
	\centering
	\includegraphics[scale=0.72]{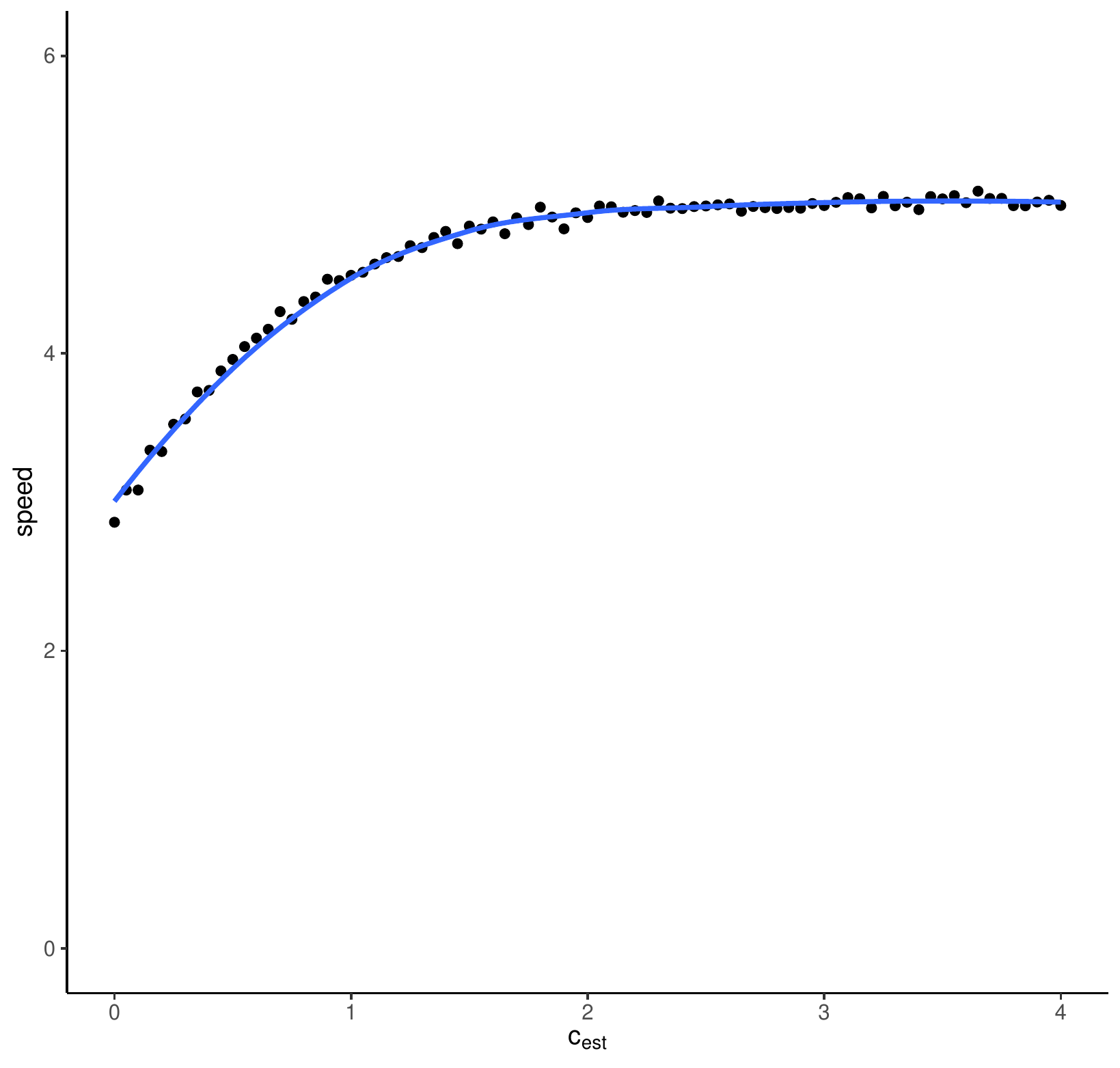}
	\caption{The speed as a function of  $ c_{\text{est}}$.
	The other parameter with is fixed $c_{\text{fec}} = 1 $.
Each estimate is computed  as $\frac{X_{t_2} - X_{t_1}}{t_2 - t_1}$
with
$t_1 = 100$ and $t_2 = 10000$.}
	\label{f2}
\end{figure}

Interestingly, we observe that for the values of $c_{\text{fec}} $
close to one, the speed increases as a function of $c_{\text{est}} $.
This phenomenon is more apparent on Figure \ref{f2},
where the speed is computed as a function of $c_{\text{est}} $
with $c_{\text{fec}} = 1 $.
This example demonstrates that 
the answer to
 Question \ref{wrong conjecture} is negative
without additional assumptions on $b_1$ and  $b_2$.

\begin{figure}[!h]
	\centering
	\includegraphics[scale=0.62]{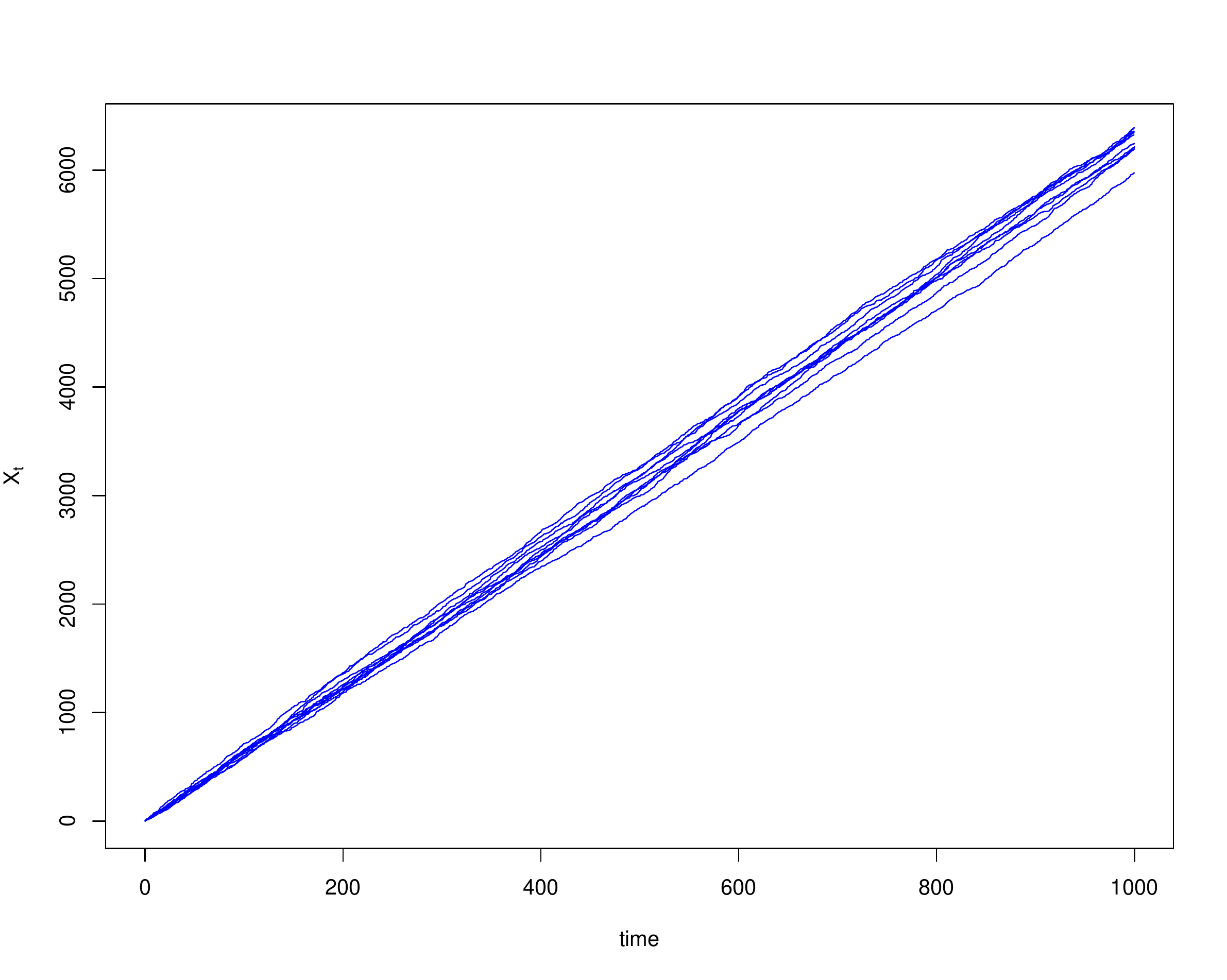}
	\caption{Ten trajectories of the tip. The parameters are $c_{\text{fec}}  =c_{\text{est}}  = 0.5$. }
	\label{f3}
\end{figure}

On Figure \ref{f3} ten different trajectories 
with 
$c_{\text{fec}}  =c_{\text{est}}  = 0.5$ are shown. 
Numerical
analysis was conducted in \textbf{\textsf{R}} \cite{R}
and figures were produced using the package ggplot2 \cite{ggplot2}.

\section*{Acknowledgements}

Viktor Bezborodov and  Tyll Krueger  are grateful for the support of the 
ZIF Cooperation group "Multiscale modelling of tumor evolution, progression and growth".
Viktor Bezborodov is grateful for the support 
from the University of Verona.
Luca di Persio thanks the "GNAMPA - Gruppo Nazionale per l'Analisi Matematica, 
la Probabilit\'a e le loro Applicazioni" 
(National Italian Group for Analysis, Probability  and their Applications)
for the received support.
Tyll Krueger is grateful for the support of the University of Verona 
during his visit. The authors would like to thank the anonymous reviewer
for their suggestions  which  have helped to improve the paper.

\bibliographystyle{alpha}
\bibliography{Sinus}

\newcommand{\etalchar}[1]{$^{#1}$}
\begin{thebibliography}{BDPK{\etalchar{+}}17}

\bibitem[ABBS13]{aidekon2013branching}
Elie A{\"\i}d{\'e}kon, Julien Berestycki, {\'E}ric Brunet, and Zhan Shi.
\newblock Branching brownian motion seen from its tip.
\newblock {\em Probability Theory and Related Fields}, 157(1-2):405--451, 2013.

\bibitem[ABK13]{arguin2013extremal}
Louis-Pierre Arguin, Anton Bovier, and Nicola Kistler.
\newblock The extremal process of branching brownian motion.
\newblock {\em Probability Theory and related fields}, 157(3-4):535--574, 2013.

\bibitem[ABR09]{Min_of_BRWs}
Louigi Addario-Berry and Bruce Reed.
\newblock Minima in branching random walks.
\newblock {\em Ann. Probab.}, 37(3):1044--1079, 2009.

\bibitem[Ahl19]{Ahl_FPP_19}
Daniel Ahlberg.
\newblock A temporal perspective on the rate of convergence in first-passage
  percolation under a moment condition.
\newblock {\em Braz. J. Probab. Stat.}, 33(2):397--401, 2019.

\bibitem[A{\"\i}d13]{aidekon2013convergence}
Elie A{\"\i}d{\'e}kon.
\newblock Convergence in law of the minimum of a branching random walk.
\newblock {\em The Annals of Probability}, 41(3A):1362--1426, 2013.

\bibitem[AMP02]{frogshape}
O.~S.~M. Alves, F.~P. Machado, and S.~Yu. Popov.
\newblock The shape theorem for the frog model.
\newblock {\em Ann. Appl. Probab.}, 12(2):533--546, 2002.

\bibitem[BDF{\etalchar{+}}]{fec_19}
Viktor {Bezborodov}, Luca {Di Persio}, Dmitri {Finkelshtein}, Yuri
  {Kondratiev}, and Oleksandr {Kutoviy}.
\newblock {Fecundity regulation in a spatial birth-and-death process}.
\newblock {\em Stoch. Dyn.}
\newblock To appear.

\bibitem[BDPK{\etalchar{+}}17]{shapenodeath}
V.~Bezborodov, L.~Di~Persio, T.~Krueger, M.~Lebid, and T.~O\.za\'nski.
\newblock Asymptotic shape and the speed of propagation of continuous-time
  continuous-space birth processes.
\newblock {\em Advances in Applied Probability}, 50(1):74–101, 2017.

\bibitem[BDPKT]{trunc_and_crop}
V.~Bezborodov, L.~Di~Persio, T.~Krueger, and P.~Tkachov.
\newblock Spatial growth processes with long range dispersion: microscopics,
  mesoscopics, and discrepancy in spread rate.
\newblock {\em Ann. Appl. Probab.}
\newblock To appear.

\bibitem[Big95]{Big95}
J.~D. Biggins.
\newblock The growth and spread of the general branching random walk.
\newblock {\em Ann. Appl. Probab.}, 5(4):1008--1024, 1995.

\bibitem[Blo13]{Blo13}
O.~Blondel.
\newblock Front progression in the {E}ast model.
\newblock {\em Stochastic Process. Appl.}, 123(9):3430--3465, 2013.

\bibitem[CG08]{CG08}
Patrick Cattiaux and Arnaud Guillin.
\newblock Deviation bounds for additive functionals of {M}arkov processes.
\newblock {\em ESAIM Probab. Stat.}, 12:12--29, 2008.

\bibitem[Che04]{ChenMarkov}
Mu-Fa Chen.
\newblock {\em From {M}arkov chains to non-equilibrium particle systems}.
\newblock World Scientific Publishing Co., Inc., River Edge, NJ, second
  edition, 2004.

\bibitem[Dei03]{Dei03}
M.~Deijfen.
\newblock Asymptotic shape in a continuum growth model.
\newblock {\em Adv. in Appl. Probab.}, 35(2):303--318, 2003.

\bibitem[Dur79]{Dur79}
Richard Durrett.
\newblock Maxima of branching random walks vs. independent random walks.
\newblock {\em Stochastic Processes and their Applications}, 9(2):117--135,
  1979.

\bibitem[Dur83]{Dur83}
R.~Durrett.
\newblock Maxima of branching random walks.
\newblock {\em Z. Wahrsch. Verw. Gebiete}, 62(2):165--170, 1983.

\bibitem[Dur88]{Dur88}
Richard Durrett.
\newblock {\em Lecture notes on particle systems and percolation}.
\newblock The Wadsworth \& Brooks/Cole Statistics/Probability Series. Wadsworth
  \& Brooks/Cole Advanced Books \& Software, Pacific Grove, CA, 1988.

\bibitem[FKK13]{est_fec}
Dmitri Finkelshtein, Yuri Kondratiev, and Oleksandr Kutoviy.
\newblock Establishment and fecundity in spatial ecological models: statistical
  approach and kinetic equations.
\newblock {\em Infin. Dimens. Anal. Quantum Probab. Relat. Top.}, 16(2):24,
  2013.

\bibitem[Gan00]{Gan00}
Nina Gantert.
\newblock The maximum of a branching random walk with semiexponential
  increments.
\newblock {\em Ann. Probab.}, 28(3):1219--1229, 2000.

\bibitem[Gar95]{Garcia}
N.~L. Garcia.
\newblock Birth and death processes as projections of higher-dimensional
  poisson processes.
\newblock {\em Adv. in Appl. Probab.}, 27(4):911–930, 1995.

\bibitem[GK06]{GarciaKurtz}
N.~L. Garcia and T.~G. Kurtz.
\newblock Spatial birth and death processes as solutions of stochastic
  equations.
\newblock {\em ALEA Lat. Am. J. Probab. Math. Stat.}, 1:281--303, 2006.

\bibitem[Hel82]{CLT_mart_CT}
Inge~S. Helland.
\newblock Central limit theorems for martingales with discrete or continuous
  time.
\newblock {\em Scand. J. Statist.}, 9(2):79--94, 1982.

\bibitem[HH80]{HH80}
P.~Hall and C.~C. Heyde.
\newblock {\em Martingale limit theory and its application}.
\newblock Academic Press, Inc. [Harcourt Brace Jovanovich, Publishers], New
  York-London, 1980.
\newblock Probability and Mathematical Statistics.

\bibitem[IW89]{IkedaWat}
N.~Ikeda and Sh. Watanabe.
\newblock {\em Stochastic differential equations and diffusion processes},
  volume~24 of {\em North-Holland Mathematical Library}.
\newblock North-Holland Publishing Co., Amsterdam; Kodansha, Ltd., Tokyo,
  second edition, 1989.

\bibitem[Kes93]{Kes_fpp_93}
Harry Kesten.
\newblock On the speed of convergence in first-passage percolation.
\newblock {\em Ann. Appl. Probab.}, 3(2):296--338, 1993.

\bibitem[KS08]{KesS08}
H.~Kesten and V.~Sidoravicius.
\newblock A shape theorem for the spread of an infection.
\newblock {\em Ann. of Math. (2)}, 167(3):701--766, 2008.

\bibitem[Lig85]{Lig85subadd}
T.~M. Liggett.
\newblock An improved subadditive ergodic theorem.
\newblock {\em Ann. Probab.}, 13(4):1279--1285, 1985.

\bibitem[LZ15]{CLT_for_CTMC}
Yuanyuan Liu and Yuhui Zhang.
\newblock Central limit theorems for ergodic continuous-time {M}arkov chains
  with applications to single birth processes.
\newblock {\em Front. Math. China}, 10(4):933--947, 2015.

\bibitem[{R C}19]{R}
{R Core Team}.
\newblock {\em R: A Language and Environment for Statistical Computing}.
\newblock R Foundation for Statistical Computing, Vienna, Austria, 2019.

\bibitem[Ric73]{Rich73}
D.~Richardson.
\newblock Random growth in a tessellation.
\newblock {\em Proc. Cambridge Philos. Soc.}, 74:515--528, 1973.

\bibitem[vdG95]{vdG95}
Sara van~de Geer.
\newblock Exponential inequalities for martingales, with application to maximum
  likelihood estimation for counting processes.
\newblock {\em Ann. Statist.}, 23(5):1779--1801, 1995.

\bibitem[Wic09]{ggplot2}
Hadley Wickham.
\newblock {\em ggplot2: Elegant Graphics for Data Analysis}.
\newblock Springer-Verlag New York, 2009.

\bibitem[Wu00]{Wu00}
Liming Wu.
\newblock A deviation inequality for non-reversible {M}arkov processes.
\newblock {\em Ann. Inst. H. Poincar\'{e} Probab. Statist.}, 36(4):435--445,
  2000.

\end{thebibliography}

\end{document}